\documentclass[11pt]{amsart}
\usepackage{amssymb}
\usepackage{graphicx}
\usepackage{amsmath, color}
\usepackage{algorithmic}
\usepackage{psfrag}
\def\h{\mathbf{h}}
\def\F{\mathbf{F}}
\def\f{\mathbf{f}}
\def\G{\mathbf{G}}
\def\ge{G_{\rightarrow}}
\def\dge{\dot{G}_{\rightarrow}}
\def\gs{G_{\downarrow}}
\def\dgs{\dot{G}_{\downarrow}}
\def\R{\mathbf{R}}
\def\rs{R_{\downarrow}}
\def\drs{\dot{R}_{\downarrow}}
\def\re{R_{\rightarrow}}
\def\dre{\dot{R}_{\rightarrow}}

\begin{document}

\title{Full and half Gilbert tessellations with rectangular cells}

\author{James Burridge, Richard Cowan, Isaac Ma }

\date{\today}

\maketitle

\noindent \begin{center}  *{\scriptsize Department of Mathematics,
University of Portsmouth,
Portsmouth, UK. \textsf{james.burridge@gmail.com} }\\
**{\scriptsize School of Mathematics and Statistics, University of
Sydney, NSW, 2006, Australia. \textsf{rcowan@usyd.edu.au} }\\
***{\scriptsize Lee Wai Lee Vocational Studies Institute, Hong Kong.
\textsf{isaacma@vtu.edu.hk} }\end{center}

\vspace{.6cm}
\begin{abstract}
We investigate the ray--length distributions for two different
rectangular versions of Gilbert's tessellation \cite{eng}. In the
\emph{full} rectangular version, lines extend either horizontally
(with east- and west--growing rays) or vertically (north- and
south--growing rays) from seed points which form a Poisson point
process, each ray stopping when another ray is met. In the
\emph{half} rectangular version, east and south growing rays do not
interact with west and north rays. For the half rectangular
tessellation we compute analytically, via recursion, a series
expansion for the ray--length distribution, whilst for the full
rectangular version we develop an accurate simulation technique,
based in part on the stopping--set theory of Zuyev \cite{zuy99},  to
accomplish the same. We demonstrate the remarkable fact that plots
of the two distributions appear to be identical when the intensity
of seeds in the half model is twice that in the full model. Our
paper explores this coincidence mindful of the fact that, for one
model, our results are from a simulation (with inherent sampling
error).  We go on to develop further analytic theory for the
half--Gilbert model using stopping--set ideas once again, with some
novel features. Using our theory, we obtain exact expressions for
the first and second moment of ray length in the half--Gilbert
model. For all practical purposes, these results can be applied to
the full--Gilbert model --- as much better approximations than those
provided by Mackissack and Miles \cite{mm}.
\end{abstract}

\section{Introduction}\label{Intro}

Consider a stationary Poisson point process in the plane, of
intensity $\lambda$. The particles of this process are called
\emph{seeds}, aptly so because at a given time $t=0$ they each
initiate the growth of a line. The directions of the lines are
randomly distributed, uniformly on $(0,\pi]$, and independent of
each other and of the seed locations. Each line grows
bidirectionally from its seed at the same rate; thus two \emph{rays}
grow from each seed. When a ray encounters a line that has already
grown across its path, the growth of that ray stops. Eventually all
rays stop growth and a tessellation of the plane is formed.

The completed structure has become known as the Gilbert tessellation
after Edgar N. Gilbert. It is notoriously difficult to analyse and
even the expected length of a typical completed ray has not been
found. There is no published paper by Gilbert on the topic; notes he
supplied appear in a book by Noble, with due acknowledgement to
Gilbert. Citations have typically attributed the notes to Gilbert
(as we do, see \cite{eng}).

A version of the model where the directions of growth were confined
to two orthogonal directions, vertical ($V$) and horizontal ($H$),
was discussed by Mackissack and Miles \cite{mm}. A tessellation of
the plane by rectangles results in their model. This structure too
has not yielded to analysis, although when seeds are equally likely
to be $V$ or $H$ the authors did provide an analytic approximation
(based on ideas of Gilbert) to the expected ray length, namely
$\sqrt {2/\lambda}$. The merits of this approximation have not been
evaluated in the literature to date.

\begin{figure}[h]
    \begin{center}
        \includegraphics[width=70mm]{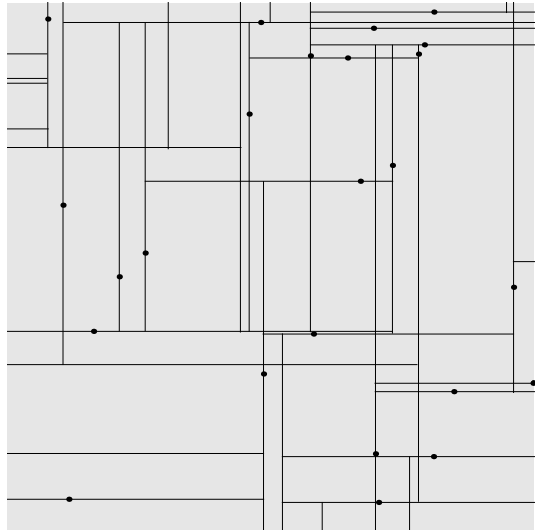}
        \caption{\scriptsize \label{halfGilbert} The Cowan--Ma (or half--Gilbert)
        rectangular tessellation when $V$--type and $H$--type seeds have
        equal intensities.   }
    \end{center}
\end{figure}

The current paper arises from work done in 1997 by the second and
third authors (Cowan and Ma). They obtained some analytic results
for an even simpler $V$\&$H$--model, whereby the growth of
eastward--growing rays is halted only by southward--growing rays
(and vice versa). Westward and northward have the same reciprocity.
A realisation of their tessellation is given in Figure
\ref{halfGilbert}.

Cowan and Ma placed a recurrence relationship (see (\ref{recur})
below) from their work on the internet \cite{cm}, though without
proof. The background to this recurrence is as follows.

\begin{figure}[h]
        \psfrag{O}{$O$}
        \psfrag{P}{$P$}
    \psfrag{A}{$A$}
        \psfrag{B}{$B$}
    \psfrag{C}{$C$}
        \psfrag{D}{$D$}
    \psfrag{E}{$E$}
    \psfrag{Q}{$Q$}
 \begin{center}
        \includegraphics[width=130mm]{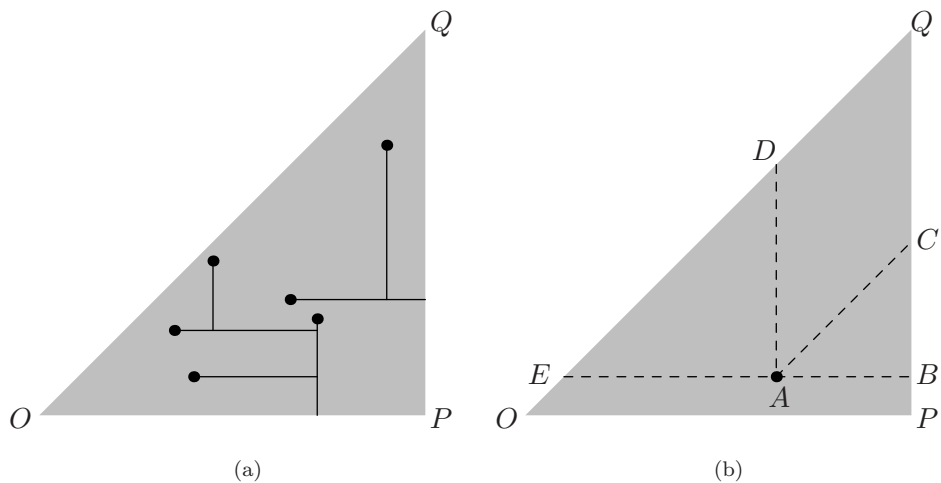}\\
                {\scriptsize \hspace{.5cm}(a)\hspace{6cm}(b)}
        \caption{\scriptsize \label{CMTri} Diagrams to assist the proof of the Cowan--Ma recurrence. }
    \end{center}
\end{figure}

Consider the isosceles right--angle triangle $POQ$ in Figure
\ref{CMTri}(a). Here $|OP|=|PQ|=\ell$. Suppose $n$ seeds lie inside
the triangle, uniformly and independently distributed; the figure
uses $n=6$. East or south growth of the rays is shown. Because of
the blocking rules, only some of the rays reach the boundary of the
triangle $POQ$.

Cowan and Ma investigated the probability $h_n$ that no rays hit the
boundary within the segment $OP$. This can also be interpreted as
the probability that $L$, the final length of a test ray commencing
eastward growth from $O$, is $> \ell$.

Their recurrence relationship for $h_n$ was as follows. For $n\geq
1$,
\begin{align}
    h_n =& \frac{n!\,q}{(2n)!}\sum_{u=0}^{n-1}\sum_{v=0}^{n-1-u}
    \frac{2^{n-u-v}\ h_u\ h_v\ (n-1+u-v)!\ (n-1-u+v)!}{u!\ v!
    \ (n-1-u-v)!}, \label{recur}
\end{align}
with $h_0=1$. Here $q$ is the proportion of seeds which grow
horizontally. The recurrence does not involve $\ell$, so $h_n$ does
not depend on $\ell$ --- as is obvious from the scale invariance of
the problem posed by Figure \ref{CMTri}(a).

This recurrence is a useful analytic step, providing precise
information on $\mathbb{E}(L)$ and $F(\ell) := \Pr\{L\leq \ell\}$.
\begin{equation}
\label{ES_dist} F(\ell) = 1 - \sum_{n\geq 0} h_n \frac{(\lambda
\ell^2)^n \exp(-\lambda \ell^2/2)}{2^n n!},
\end{equation}
from which we deduce (in an extended notation which includes
$\lambda$) that $F_{\lambda}(\ell) = F_1(\sqrt{\lambda}\,\ell)$.
Also
\begin{align*}
\mathbb{E}(L) = & \int_0^\infty [1-F(\ell)] d\ell \\
 = & \sum_{n \geq 0} \frac{h_n}{ n!} \int_0^\infty (\lambda \ell^2/2)^n
   \exp(-\lambda \ell^2/2) d\ell \\
 = & \frac{1}{\sqrt{2\lambda}} \sum_{n \geq 0} \frac{h_n
 \Gamma(n+\frac{1}{2})}{n!}.
\end{align*}

In Section \ref{CMrecur}, we report the proof used to derive the
recurrence relationship (\ref{recur}) and plot the probability
density function of the random variable $L$. The plot has an
extraordinary property, discovered when certain simulations of the
full rectangular--Gilbert model done by our first author, Burridge
\cite{bur}, were also plotted. The probability density function of
the Cowan--Ma model with $\lambda = 2$ was indiscernible from that
of the full rectangular--Gilbert model with $\lambda = 1$.

Section \ref{DJBSim} presents Burridge's simulation study, that has
a very high level of accuracy, and discusses this surprise
coincidence
--- which raises somewhat the profile of the Cowan--Ma model. As
well as having interest in its own right as a tessellation model
with tractable mathematics, the model provides approximations for
the full--Gilbert rectangular model. For example, the Cowan--Ma
model
--- which we also called the half--Gilbert model because it has half
of the blocking mechanisms --- provides a much better approximation
for $\mathbb{E}(L)$ in the full model when compared with the
Mackissack/Miles approximation, which is $\mathbb{E}(L)\approx \sqrt
{2/\lambda}$  when $q=\tfrac12$.

In Section \ref{SSHalf} our work pushes further the tractability of
the half--Gilbert model finding; most notably, we find that the mean
ray length when $q=\tfrac12$ is given by the formula:
\begin{equation*}
 \mathbb{E}(L)= \frac{\pi}{\sqrt{\lambda}
\left(\Gamma(\tfrac{3}{4})\right)^2}.
\end{equation*}

In both our simulation and analytical work we have employed Zuyev's concept of stopping set sequences \cite{zuy99} and the distributional results
for the areas of these sets. To achieve the analytic results, we have incorporated a
new concept into the analysis, the idea of dead zones which
influence the formation of the next stopping set in the sequence. Our most complete analysis is for the balanced case, $q=\tfrac12$,
because some results become rather complicated when $q \ne
\tfrac12$. The expected ray length in the latter case is reported, without proof, in the appendix.

\section{The Cowan--Ma recurrence relation}\label{CMrecur}

We now prove (\ref{recur}) for general $q$. Obviously, $h_0=1$ and
$h_1=\frac{1}{2}$. When $n\geq1$, we label the seed closest to $OP$
as $A$. See Figure 2(b). If the distance from $A$ to $OP$ is denoted
by the random variable $Y$, it is easily shown that $Y$ has
probability density function $g_Y(y)= 2n(1-y)^{2n-1},\  0\leq y\leq
1$. Furthermore, the conditional probability density function of
$X:=|AB|$ given $Y$ is
\begin{align*}
    g(x|y) =&\ \frac{1}{1-y},\qquad 0\leq x\leq 1-y.
\end{align*}
Denote the event that no rays hit $OP$ by $\mathcal{E}_n$. Then
\begin{align*}
    \Pr\{\mathcal{E}_n | x,y\} =&\ \Pr\{\mathrm{seed\ }A\mathrm{\ grows
    \ eastward\ and\ reaches\ }B\\
    &\qquad \mathrm{\ and\ no\ ray
    \ grows\ across\ the\ segment\ }EA\} \\
    =&\ q \Pr\{\mathrm{no\ ray\ grows\ across\ }AB
    \mathrm{\ \emph{and}\ no\ ray
    \ grows\ across\ }EA\}.
\end{align*}
To evaluate the right--hand side, we partition the domain above $EB$
into the three zones that are shown in Figure 2(b). We then consider
the trinomial distribution by which the remaining $(n-1)$ seeds are
allocated to these zones: $u$ to $ABC$, $v$ to $EAD$ and the
remaining $(n-1-u-v)$ to $ACQD$. This leads, for each $(u,v)$, to a
rather pleasing representation of the problem into two problems
self--similar to the original one. Continuing, using $||\cdot||$ as
area, we write $\Pr\{\mathcal{E}_n | x,y\}$ as
\begin{align*}
    &\ q \sum_{u=0}^{n-1}\sum_{v=0}^{n-1-u}
    \frac{(n-1)!\ ||ABC||^u\ \ ||EAD||^v
    \ \ ||ACQD||^{n-1-u-v}}{u!\ v!\ (n-1-u-v)!\ \ ||EBQ||^{n-1}}\ \ \times\\
    &\qquad\Pr\{\mathrm{no\ ray\ grows\ across\ }AB
    \mathrm{\ \emph{and}\ no\ ray
    \ grows\ across\ }EA|u,v\}\\
    =&\ q \sum_{u=0}^{n-1}\sum_{v=0}^{n-1-u}
    \frac{(n-1)!\ (\frac{x^2}{2})^u\ \ \bigl(\frac{(1-x-y)^2}{2}\bigr)^v
    \ \ [x(1-x-y)]^{n-1-u-v}}{u!\ v!\ (n-1-u-v)!\ \ (\frac{(1-y)^2}{2})^{n-1}}
    \ \ \times\\
    &\qquad\Pr\{\mathrm{no\ ray\ grows\ across\ }AB|u\}
    \Pr\{\mathrm{no\ ray
    \ grows\ across\ }EA|v\}\\
    =&\ q \sum_{u=0}^{n-1}\sum_{v=0}^{n-1-u}
    \frac{(n-1)!\ x^{2u}(1-x-y)^{2v}
    [2x(1-x-y)]^{n-1-u-v}}{u!\ v!\ (n-1-u-v)!
    \ \ (1-y)^{2(n-1)}}\ h_u\ h_v.
\end{align*}
Unconditional on $x$ and $y$, and with $n\geq 1$,
\begin{align*}
    h_n=&\ \Pr\{\mathcal{E}_n\}\\
    =&\ \int_0^1\int_0^{1-y}\ \Pr\{\mathcal{E}|x,y\}g_Y(y)g(x|y)\ dxdy\\
    =&\ q \sum_{u=0}^{n-1}\sum_{v=0}^{n-1-u}
    \frac{(n-1)!\ h_u\ h_v}{u!\ v!\ (n-1-u-v)!}\ \ \times\\
    &\qquad\quad 2n \int_0^1\int_0^{1-y}\ x^{2u}(1-x-y)^{2v}
    [2x(1-x-y)]^{n-1-u-v}
    \ dxdy\\
    =&\ (n-1)!\ q \sum_{u=0}^{n-1}\sum_{v=0}^{n-1-u}
    \frac{2^{n-1-u-v}\ h_u\ h_v}{u!\ v!\ (n-1-u-v)!}\ \ \times\\
    &\qquad \qquad 2n \int_0^1\int_0^{1-y}\ x^{n-1+u-v}(1-y-x)^{n-1-u+v}
    \ dxdy\\
    =&\  n!\ q\sum_{u=0}^{n-1}\sum_{v=0}^{n-1-u}
    \frac{2^{n-u-v}\ h_u\ h_v}{u!\ v!\ (n-1-u-v)!}\ \times\\
    &\qquad \qquad \int_0^1\ (1-y)^{2n-1}B(n+u-v,n-u+v)
    \ dy\\
    =&\ \frac{n!\ q}{(2n)!}\sum_{u=0}^{n-1}\sum_{v=0}^{n-1-u}
    \frac{2^{n-u-v}\ h_u\ h_v\ (n-1+u-v)!\ (n-1-u+v)!}{u!\ v!\ (n-1-u-v)!}.
\end{align*}
We augment this recurrence with the result $h_0=1$. This completes
the proof of (\ref{recur}). We note that the sequence $h_0, h_1,
h_2, ...$ commences $1, \tfrac12, \tfrac13, \tfrac{29}{120},
\tfrac{11}{60}, ...$ when $q=\tfrac12$.

\begin{figure}
    \psfrag{f}{$f(\ell)$}
    \psfrag{l}{$\ell$}
     \psfrag{A}{$q=\tfrac{1}{4}$}
    \psfrag{B}{$q=\tfrac{1}{2}$}
    \psfrag{C}{$q=\tfrac{3}{4}$}
\begin{center}
\includegraphics[width=11 cm]{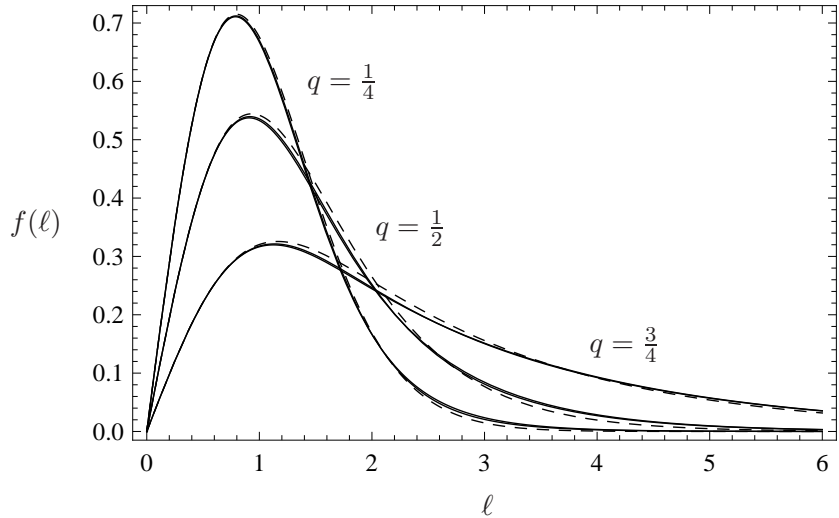}
\caption{\scriptsize \label{pdfs_with_q} The three solid curves are
the probability density functions $f$ for the final length of a
typical eastward--growing ray in the Cowan--Ma model. Each is based
on $\lambda=2$ and with three $q$ values: $\tfrac14, \tfrac12$ and
$\tfrac34$. We show later in the paper that: (a) each of these solid
curves actually comprises two curves overlaid, the second being the
curve from the full--Gilbert model, with $\lambda = 1$; (b) the
dashed lines are the probability density functions from Gilbert's
heuristic `mean field' analysis, also valid for both models.}
\end{center}
\end{figure}

The recurrence together with (\ref{ES_dist}) can be used to plot
$f(\ell):= F'(\ell)$ against $\ell$ for various values of $q$ (see
Figure \ref{pdfs_with_q}).

\section{Simulation of the full rectangular--Gilbert tessellation}
\label{DJBSim}

Finding coefficients analogous to $h_n$ for the full rectangular
model is a formidable task because of the complexity of the blocking
effects. Lacking self--similar zones akin to those discovered by
Cowan and Ma in their model, we have devised an efficient way of
accounting for these effects by simulation.

\begin{figure}[h]
   \psfrag{A}{$A$}\psfrag{B}{$B$}\psfrag{C}{$C$}
\begin{center}
\includegraphics[width=80mm]{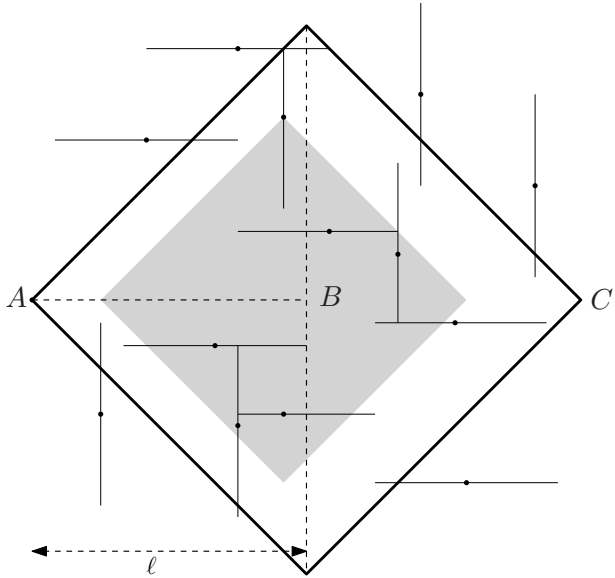}\hspace{1cm}\\
\caption{\scriptsize \label{active_square} Only particles within the
large square (whose diagonal is $AC$) can influence the event that
an $H$--type ray starting at $A$ does not reach $B$ (due to
intersection of the line--segment $AB$ by vertical rays). The role
of the smaller shaded square is described in the text.}
\end{center}
\end{figure}

The analogue of the isosceles triangle used in Figure \ref{CMTri}
 is a square, rotated so that its
diagonal $AC$ lies east--west, as illustrated in Figure
\ref{active_square}. To study the growth of horizontal rays, we
consider an $H$--type \emph{test seed} located at the western corner
of the square, marked A in the figure, and define:
\begin{equation*}
\h_n =\ \Pr\{\mathrm{ray\ from\ test\ seed\ } A \mathrm{\ reaches\ }
B \mid n \mathrm{\ seeds\ in\ the\ square}  \}.
\end{equation*}
The only seeds that can block the test ray lie in the western side
of the square, but whether or not they do so depends also on the
configuration of seeds in the eastern side. Seeds outside the square
have no influence.

By analogy with equation (\ref{ES_dist}) the ray
length distribution for the rectangular Gilbert tessellation is:
\begin{equation*}
 \F(\ell) = 1- \sum_{n\geq 0} \h_n \frac{(2 \lambda
\ell^2)^n \exp(-2\lambda \ell^2)}{n!},
\end{equation*}
from which $\F_{\lambda}(\ell) = \F_1(\sqrt{2\lambda}\,\ell)$ is
deduced.

\vspace{.4cm}\textbf{An obvious method:} The naive approach to
estimating $\h_n$ would be to repeatedly populate the large square
in Figure \ref{active_square} with $n$ seeds (each independently of
$H$--type with probability $q$) and, each time, determine if the
line $AB$ is intersected. This can be accomplished using the
following recursive algorithm which decides if a ray, extending in
compass--direction $u\in \{\rightarrow, \uparrow, \leftarrow,
\downarrow\}$ from one seed $s^*$  will be blocked within a distance
$d$. The algorithm, \textsf{block}, outputs a logical value:

\begin{equation*}
\mathsf{block}(s^*, d, u) = \left\{
\begin{array}{l l}
  \texttt{true} & \textrm{if ray is blocked}   \\
  \texttt{false} &  \textrm{if ray is not blocked}. \\
\end{array} \right.
\end{equation*}

Whether or not $s^*$ is blocked within a distance $d$ depends only
on the configuration of seeds within a square of diagonal $2d$ along
which its produced ray travels. Let the compass--direction of this
ray be $u$, and let us denote by $\Delta(s^*, d, u)$ the isosceles
triangle which forms the half of the square closest to $s^*$. Let
the type ($H$ or $V$) of seed $s$ be $t(s)$. The algorithm
\textsf{block}$(s^*, d, u)$ runs as follows:

\begin{algorithmic}
\FORALL{$s \in \Delta(s^*, d, u)$}
    \IF{$t(s) \neq t(s^*)$}
        \STATE compute the perpendicular distance, $d_s$, and compass--direction,
        $u_s$ from $s$ to $s^*$'s produced ray.
           \IF{$\mathsf{block}(s,d_s,u_s) =$ false}
                \RETURN \TRUE
           \ENDIF
    \ENDIF
\ENDFOR
\RETURN \FALSE
\end{algorithmic}

For example, if $s^*$ is the $H$--type seed at $A$, then the
computer programme calls \textsf{block}$(s^*, \ell, \rightarrow)$.
This invokes recursive calls to \textsf{block} for every $V$--type
seed in the left isosceles triangle (until a \texttt{true} value is
returned by the call). In Figure \ref{active_square}, the shaded
region  with a $V$--type seed $s$ at the top shows a square that is
investigated by one of the recursive calls, specifically by the call
\textsf{block}$(s, d_s, \downarrow)$, where $2d_s$ is the diagonal
length of the shaded square.

In principle, we can conduct this simulation for each $n$ up to
(say) $300$. For each $n$, we would generate the seeds in the square
(with diagonal $AC$) $N$ times, where $N$ would be very large. An
estimate of $\h_n,\ 0\leq n\leq 300$, is thereby generated for
$H$--type rays. Then, if $q\ne \tfrac12$, we would repeat the whole
procedure for $V$--type rays. It is a lengthy process, despite the
potential saving if an early--tested seed $s$ returns \texttt{true}
--- implying that others don't have to be tested.
\begin{figure}[h]
\begin{center}
\includegraphics[width=70mm]{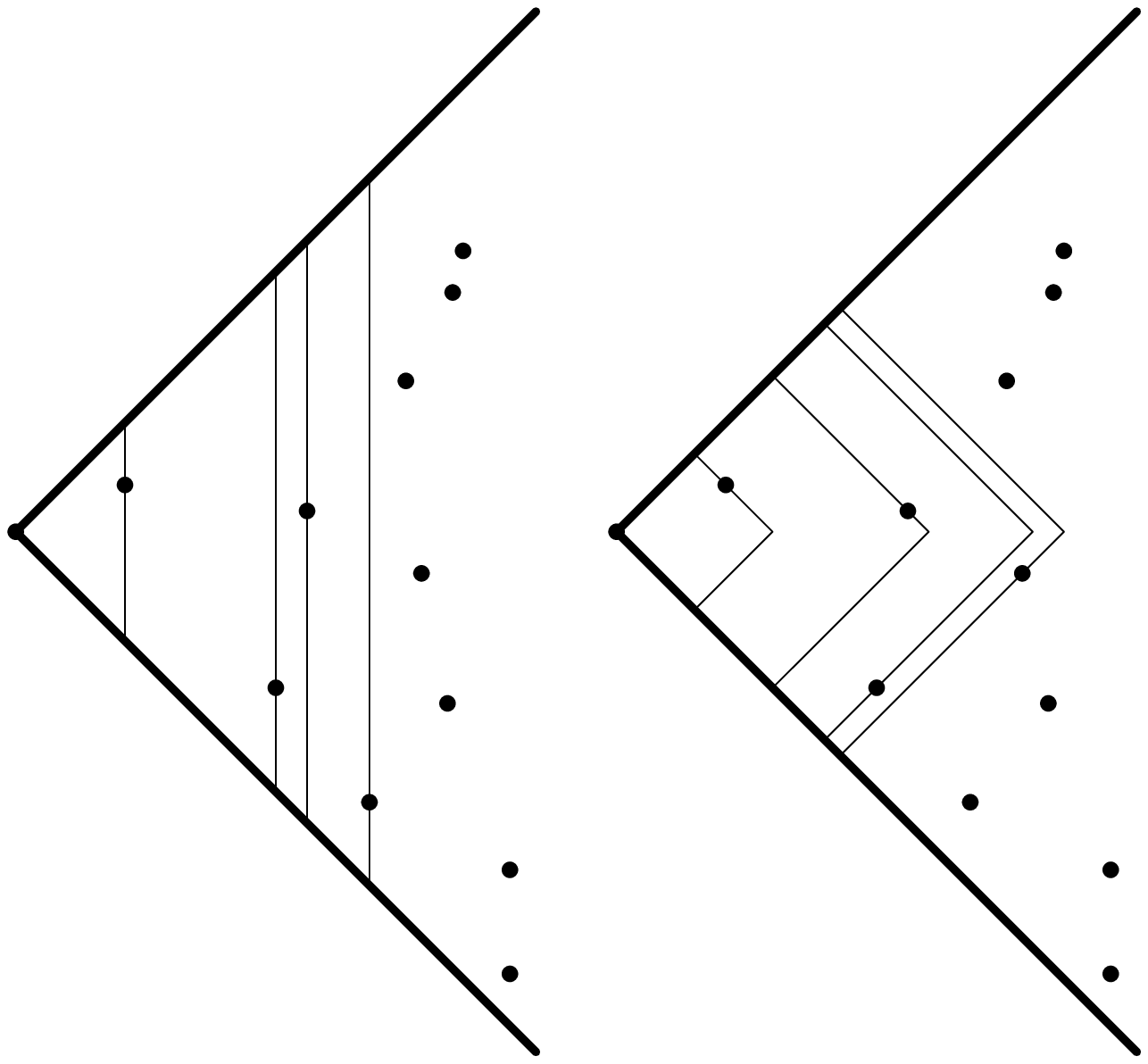}\hspace{1.5cm}
\includegraphics[width=55mm]{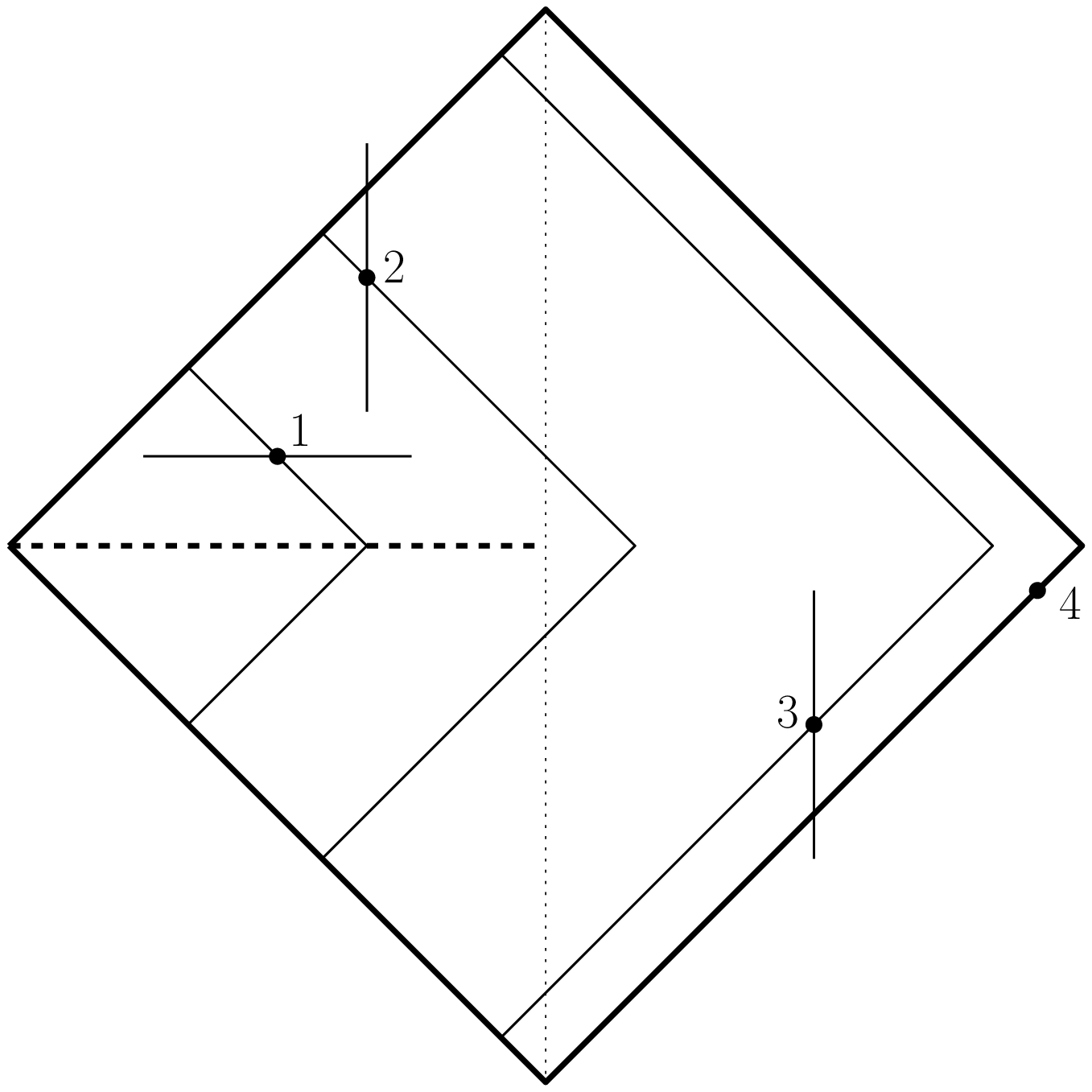}
        {\scriptsize \hspace{4.5cm}(a)\hspace{3cm}(b)\hspace{4.3cm}(c)}
\caption{\scriptsize \label{sets} Nested stopping sets are shown in
(a) and (b). In (c), the `efficient algorithm' is at step 3 and at
least one more step will be needed before we see a ray in step $n$
crossing the dashed half--diagonal of $S_{n+1}$.}
\end{center}
\end{figure}

\vspace{.4cm}\textbf{Stopping sets:} To shorten the task, we have
devised a method based on stopping sets (a concept defined by Zuyev
\cite{zuy99} and amplified in \cite{cowan03}). Consider the
unbounded quadrant that lies between the half--lines $y=x$ and
$y=-x$, with $x\geq0$, partly shown in Figure \ref{sets}(a). A
stationary Poisson process of seeds with intensity $\lambda$ exists
in the quadrant. A triangular set whose eastern boundary is vertical
and western vertex is the quadrant's apex is gradually expanded,
stopping briefly whenever its boundary hits a seed --- before
continuing its expansion. The set stopped by the $k$th seed
encountered is called $S_k$. This process creates a nested sequence
of random sets. We denote the area of $S_1$ by $E_1$ and the areas
of the region $S_k \setminus S_{k-1}$ by $E_k,\ k>1$. Another
nesting arrangement is shown in Figure \ref{sets}(b), this time with
squares  and a different ordering of the seeds.

Zuyev showed, among other things of a more general nature, that any
expanding domain constructing a nest of compact sets in the manner
described above --- through a sequence of \emph{stops} caused by
seed hits --- creates areas $E_1, E_2, E_3, ...$ which are
independent and distributed exponentially with parameter $\lambda$.
The domain might have a complicated geometry because the expansion
rule is allowed to depend on the seeds that it contains (and, being
closed, this includes seeds on the domain's boundary). In the two
examples of Figure \ref{sets}, the expansion rule is straightforward
and doesn't depend on the internal seeds.

Most importantly for the validity of Zuyev's distributional results,
 neither the expansion rule nor the stopping rule for $S_1$ should depend on
seeds \emph{outside} the expanding domain. This prohibition plays
two roles.
\begin{itemize}
    \item it helps establish that $E_1$ is exponentially distributed;
    \item it also
allows one to say that the point process of seeds outside the
stopping set $S_1$ is still a stationary Poisson point process with
unchanged intensity \emph{given} the information within $S_1$ (a
notion formalised by Theorem 2 of \cite{cowan03}).
\end{itemize}
This allows the argument to be extended sequentially to $E_2, E_3,
...$ and $S_2, S_3, ....$.

We also note that Zuyev's results are not guaranteed if
randomisations apart from the Poisson process of seeds affect the
growth and stopping. No such complication occurs in this section of
our paper, although we must address the issue later in Section
\ref{SSHalf}.

Stopping sets constructed in this way have other properties. The
$i$th seed $s_i$ is uniformly distributed on the growth frontier of
$S_i$ and the seeds $s_1, s_2, ..., s_n$ are uniformly and
independently distributed in the set $S_{n+1}$. Those of Figure
\ref{sets}(b) have a property that no other nesting has: if $s_i$ is
$V$--type, then whether or not it reaches the east--west diagonal
depends \emph{only} on seeds $s_1, s_2, ..., s_{i-1}$.

\vspace{.4cm}\textbf{Efficient algorithm:} In the context of Figure
\ref{sets}(b) with its nesting of squares, the latter property says
that the ray growth just within $S_{n+1}$ from the seeds $s_1, s_2,
..., s_n$ provides a sample of the problem that interests us ---
giving a \emph{true} or \emph{false} datum on whether a test ray is
blocked before it traverses across half the diagonal of $S_{n+1}$.
(See the illustration for $n=3$ in Figure \ref{sets}(c).) This datum
contributes to the estimation of $\h_n$. Importantly, as we show
below, if the datum is \emph{true}, then we can add a \emph{true}
datum for the estimation of all $\h_j, j>n$ --- without further
computational effort.

We start with the unbounded quadrant empty of seeds, then place an
$H$--type test seed at the apex of the quadrant. We generate the
exponentially distributed areas $E_1$ and $E_2$ and so construct the
squares $S_1$ and $S_2$ expanding from the apex. We randomly select
(uniformly) a seed point $s_1$ on the growth frontier (eastern
sides) of the inner square, $S_1$. Because of the properties
discussed above, this is equivalent to choosing the point uniformly
within the outer square $S_2$. If this seed grows a vertical ray
that intersects the diagonal, let the distance of the intersection
point from the apex be $X_1$. If not, set $X_1 =\infty$.

Let $A_i$ denote $||S_i||$, the area of $S_i$, and $\mathcal{E}_n$
denote the event that the line from the test seed reaches the centre
of a square populated with $n$ uniformly distributed seeds.
Obviously $\h_n = \Pr\{\mathcal{E}_n\}$.

If $X_1 < (A_2/2)^{1/2}$ then the simulation ends. There is no need
to generate more nested squares in order to simulate the events
$\mathcal{E}_n$, $n>1$ because we know that the half diagonal of
every subsequent square will be crossed at $X_1<(A_n/2)^{1/2}$.
Seeds on the boundaries of subsequent squares cannot influence this.
If the first seed does not cross the diagonal, or crosses such that
$X_1>(A_2/2)^{1/2}$, then we draw $S_3$ and pick a point $s_2$ on
the boundary of the second nested square $S_2$. We check if $s_2$
intersects the diagonal, accounting for any possible blocking
effects from $s_1$ by using the algorithm \textsf{block}. If so, we
let the distance from the
 apex to the \emph{closest}
intersection point be $X_2$, which will be $\leq X_1$. If $X_2 <
(A_3/2)^{1/2}$ then the simulation ends. If not, we add another
square $S_4$ and seed $s_3$ --- reaching the situation in Figure
\ref{sets}(c) --- and so on. We keep repeating the process ---
adding another seed and using \textsf{block} on that seed --- until
\textsf{block} indicates that the latest half--diagonal has been
hit. We then record that the event $\mathcal{E}_n$ fails to occur
for this and all higher values of $n$. The entity $\h_n$ for
eastward growing rays is the fraction of times that $\mathcal{E}_n$
occurs over many simulations. If $q\ne\tfrac12$, the complete
protocol is repeated with $q$ replaced by $(1-q)$ to give results
for southward growing rays.

To estimate the $\h_n$, $N=10^{9}$  simulations were performed,
requiring a running time of approximately one hour on a modern PC.
When $q=\tfrac12$, the largest number of nested squares created
before the simulation terminated was 917, which occurred once, and
the second largest number was 727, which also occurred once. The
mean number of squares created before termination was $5.25$. In the
($q=\tfrac12$) case, the estimate of expected length of each line
produced from a seed was:
\begin{equation}\label{estEL}
\mathbb{E}(L) = 1.467535\ (0.000029)
\end{equation}
where the bracketed number is the standard error, calculated with
due regard to the positive covariance between our estimators of
$\h_n$ and $\h_{n+k},\ k>0$.

{\sc Remark 1:} \emph{Our accurate estimate of the $\h_n$ values
allows the probability density function of the ray length to be
calculated. Because the two ray lengths coming from a particular
seed are independent, the standard convolution method leads to an
estimated distribution of the total line length arising from a
typical seed. Mackisack and Miles \cite{mm} claim that these two
ray--lengths are not independent, but we disagree.}

\vspace{.4cm} \textbf{The coincidence:} We found a remarkable
similarity between the probability density functions in the half
rectangular Gilbert model and the full rectangular Gilbert model
when the intensity of seeds in the former case was twice that of the
latter case. Indeed the plots were almost indistinguishable, so
Figure \ref{pdfs_with_q} effectively displays both $f$ and $\f$ for
various $q$, with $\lambda = 2$ or $\lambda = 1$ respectively.

We are mindful that the simulated results have sampling error,
albeit small. So we asked the question: are the two distributions
$F$ and $\F$ mathematically equal --- or just approximately so? To
answer this in the $(q=\tfrac12)$ case, we performed some rather
tedious exact calculations (details omitted)  which yielded:
\begin{equation*}
    \h_0= 1;\qquad\h_1 = \tfrac34;\qquad \h_2=\tfrac{7}{12};\qquad \h_3=
    \tfrac{7}{15}.
\end{equation*}
We then expanded both $F$ and $\F$ as Taylor series about the
origin.
\begin{align*}
    F(\ell) =&\ h_0+(h_1-h_0) \ell^2 + \tfrac12(h_0+h_2-2h_1)\ell^4 \\
    &\ \qquad \qquad+ \tfrac16(3h_1-3h_2-h_0+h_3)\ell^6 + o(\ell^7)\\
    =&\ 1-\tfrac12 \ell^2 +\tfrac16\ell^4 - \tfrac{31}{720}\ell^6
    +o(\ell^7).\\
    \F(\ell) =& \h_0 +2(\h_1-\h_0)\ell^2 + 2(\h_0-4\h_1+\h_2)\ell^4\\
    =&\ \qquad \qquad+\tfrac43(3\h_1-3\h_2-\h_0+\h_3)\ell^6 +o(\ell^7)\\
    =&\ 1-\tfrac12 \ell^2 +\tfrac16\ell^4 - \tfrac{32}{720}\ell^6
    +o(\ell^7).
\end{align*}
We see that these \emph{exact} series differ slightly in the fourth
term, so $F$ and $\F$ are not mathematical equal.

\vspace{.4cm}\textbf{``Mean field'' analysis when $q=\tfrac12$: }
Gilbert's original ``mean field'' analysis, which was adapted by
Mackisack and Miles \cite{mm} to the $(q=\tfrac12)$  rectangular
case, involved the rough approximation that ray ends (there being
two per seed) were uniformly spread across the plane. With this
assumption, it was possible to approximate at time $t$ the expected
number of ray ends lying within a small distance $\delta x$ of rays
that would block the growth of these ends within the next $\delta
t$.

Mackisack and Miles analyzed the $(q=\tfrac12)$ full model using two
quantities; $\R(t)$, the expected total length of rays per unit
area; $\G(t)$, the expected number of growing ends per unit area.
Recounting their work, these quantities are related exactly by
$\dot{\R}=\G$, assuming unit growth rate, and heuristically in the
full rectangular case by $\dot{\G} \approx -\tfrac{1}{2}\R\G$, with
initial conditions $\R(0)=0$ and $\G(0)=2 \lambda$. Solving these
differential equations, they found that $\G(t)\approx 2 \lambda
\textrm{sech}^2 \sqrt{\tfrac{\lambda}{2}} t$. If $L$ is the final
length of a test ray in their full Gilbert model, then:
\begin{equation}\label{MMX}
\Pr(L>\ell) =  \frac{\G(\ell)}{\G(0)}= \frac{\G(\ell)}{2
\lambda}\approx \textrm{sech}^2 \sqrt{\tfrac{\lambda}{2}} t.
\end{equation}
The expected $L$ when $q=\tfrac12$ is therefore approximated by
$\sqrt{2/\lambda}= 1.41421$ at $\lambda = 1$. This is not especially
close to the value shown in (\ref{estEL}). The solution for $\R$ was
$\R(t)\approx 2\sqrt{2\lambda}\tanh(\sqrt{2/\lambda}\,t),\ t>0$.

We have modified the analysis in \cite{mm} to deal with the
$(q=\tfrac12)$ Cowan--Ma model. We put $\dot{G}\approx-\tfrac{1}{4}R
G$ since each of the four directions of growing lines can only be
blocked by one other line type. Solving the new equation pair, we
find that the number of growing lines per unit area at time $t$ for
the half model is : $G(t) \approx 2 \lambda \textrm{sech}^2
\sqrt{\tfrac{\lambda}{4}} t$. Also $R(t)\approx
2\sqrt{2\lambda}\tanh(\sqrt{2/\lambda}\,t),\ t>0$. Furthermore
(\ref{MMX}) still holds, with $G$ replacing $\G$. So, setting
$\lambda=2$ in the half system and $\lambda=1$ in the full system we
obtain identical approximations to the probability density function
for ray length:
\begin{equation*}
 f_{2}(\ell) \approx \sqrt{2}\ \textrm{sech}^2
\frac{\ell}{\sqrt{2}}\, \tanh \frac{\ell}{\sqrt{2}}\approx
\f_{1}(\ell).
\end{equation*}
The expected ray length is: $\mathbb{E}(L)\approx \sqrt{2}$.

So we have shown that the mean field approximations in the two
models are equal, when $q=\tfrac12$. Indeed, our analysis for
$q\ne\tfrac12$, developed in the next sub--section, shows that the
two approximations are also equal when $q\ne \tfrac12$.

\vspace{.4cm}\textbf{Mean field analysis when $q\ne \tfrac12$:} When
the intensities of $H$- and $V$--type seeds are not equal, the rays
of east--growing and south--growing have different length
distributions. So a system of four differential equations and four
initial values is needed, in variables (for the half--Gilbert model)
$\gs, \ge, \rs$ and $\re$.
\begin{alignat*}{2}
    \dre(t)=&\ \ge(t)&\qquad\dge(t)&\approx - \rs(t)\ge(t)\\
    \drs(t)=&\ \gs(t)&\qquad\dgs(t)&\approx - \re(t)\gs(t),
\end{alignat*}
combined with:
\begin{align*}
    \re(0) =&\ \rs(0) = 0;\qquad\qquad
    \ge(0) =\ q\lambda;\qquad \gs(0)=\ (1-q)\lambda.
\end{align*}
Replacing $\approx$ with $=$ and eliminating $G_H$ and $G_V$, the
differential equations become
\begin{align*}
    \ddot{R}_{\rightarrow}(t)&= - \rs(t)\dre(t)\\
    \ddot{R}_{\downarrow}(t)&= -\re(t)\drs(t),
\end{align*}
augmented by
\begin{align*}
    \re(0) =&\ \rs(0) = 0\qquad\qquad
    \dre(0) =\ q\lambda\qquad \drs(0)=\ (1-q)\lambda.
\end{align*}

 We  have only
been able to solve this coupled system in series form and, even
then, with no general term recognised. Using the abbreviations
$Q:=q\lambda$ and $P:=(1-q)\lambda$,
\begin{align}\label{seriesMF}
    \re(t)=&\
    \frac{Q}{1!}\,t-\frac{PQ}{3!}\,t^3+\frac{PQ(3P+Q)}{5!}\,t^5
    -\frac{PQ(15P^2+16PQ+3Q^2)}{7!}\,t^7 \notag \\
    &\quad +\frac{PQ(105P^3+241P^2Q+135PQ^2+15Q^3)}{9!}\,t^9
    -... ,
\end{align}
with $\ge(t)$ being $\dre(t)$ (easily calculated from
\ref{seriesMF}). A \emph{Mathematica} routine to compute as many
terms as required is available from the authors. For $\rs$ and
$\gs$, simply interchange $P$ and $Q$. Note that west--growing rays
have results identical to east--growing --- likewise north and south
results are identical.

For the full--Gilbert model, the equations are very similar, but
cast in terms of the four variates $\G_V, \G_H, \R_V$ and $\R_H$.
\begin{alignat*}{2}
    \dot{\R}_H(t)=&\ \G_H(t)&\qquad\dot{\G}_H(t)&\approx - \R_V(t)\G_H(t)\\
    \dot{\R}_V(t)=&\ \G_V(t)&\qquad\dot{\G}_V(t)&\approx - \R_H(t)\G_V(t),
\end{alignat*}
combined with
\begin{align*}
    \R_H(0) =&\ \R_V(0) = 0\qquad\qquad
    \G_H(0) =\ 2q\lambda\qquad \G_V(0)=\ 2(1-q)\lambda.
\end{align*}
This leads to a solution for $\R_H(t)$ equal to the right--hand side
of (\ref{seriesMF}), but with $Q=2q\lambda$ and $P=2(1-q)\lambda$.
Thus it becomes obvious that $\R_H(t)$ with $\lambda=1$ equals
$\re(t)$ with $\lambda=2$. Likewise for the other linked pairs of
variables! Therefore, when $q\ne \tfrac12$, the two ray length
distributions (for $H$ and $V$ rays) for the full model having
intensity $\lambda$ are equal to the corresponding ray length
distributions for the half--Gilbert model with seed--intensity
$2\lambda$. All of these entities are, of course, only approximate
solutions to the true Gilbert models.

Figure \ref{pdfs_with_q} shows that their value as approximations
for the full--Gilbert model is quite good, but not nearly as good as
the analytic answers adopted from the half--Gilbert model. In the
last section of the paper, we provide more of these \emph{answers},
demonstrating that the half--Gilbert model of Cowan and Ma is
encouragingly tractable.

\section{Stopping sets and dead zones in the half--Gilbert model.}
\label{SSHalf}

It is possible to use the stopping--set concept to find exact
expressions for the first, second and in principle higher moments of
the ray length in the Cowan--Ma model. The balanced case,
$q=\tfrac12$, is easier to describe --- and that is now our focus.
We give some results for the general case in an Appendix.

\vspace{.4cm}\textbf{A different construction of stopping sets:}
Suppose that a stationary Poisson process of intensity $\lambda$
exists in the plane, with seeds marked either $H$ (east growing) or
$V$ (south growing) with equal probability. In Section \ref{DJBSim}
we have described how a nest of Zuyev's stopping sets is created
when the growth frontier of an expanding domain hits the seeds. For
the Cowan--Ma model, the seeds that are relevant for an
east--growing test ray commencing at $O$ in Figure \ref{CMTri}(a) is
the shaded region in that figure
--- or, more precisely, the \emph{unbounded} octant lying between
$y=x$ and
$y=0$, with $x\geq 0$: we call this region, the initial \emph{live
zone}.

\begin{figure}[h]
    \psfrag{r}{$r$}\psfrag{A}{$S_1$}\psfrag{B}{$D_1$}\psfrag{C}{$S_2$}
        \psfrag{y}{$y$}\psfrag{P}{$O\ $}
\begin{center}
\includegraphics[width=8 cm]{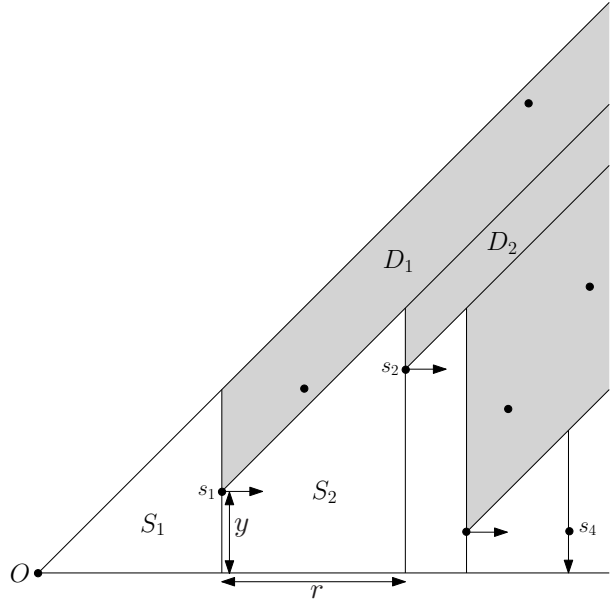}
\caption{\scriptsize \label{trapSS} Trapezoidal stopping sets and
dead zones in the half model. }
\end{center}
\end{figure}

As before we start by expanding a domain --- an isosceles right
angle triangle in this case (see Figure \ref{trapSS}) --- into the
live zone, stopping when it hits the first seed $s_1$ whose
coordinates relative to $O$ are $(x_1, y_1)$. This creates a domain
$S_1$ with area $E_1$ that is exponentially distributed. If $s_1$ is
$V$--type, then it will provide the ray that blocks the test seed;
thus $L=x_1$ and no other seeds need be considered.

Alternatively if $s_1$ is $H$--type, then, instead of growing $S_1$
(retaining its shape as an isosceles right--angle triangle and
constructing the familiar Zuyev \emph{nest} of stopping sets), we
introduce a significant modification. We remove a part of the live
zone: a `dead zone' labelled $D_1$ (see Figure \ref{trapSS}) which
has now become irrelevant, as we shall soon see.

As $S_1 \cup D_1$ has been constructed without drawing upon any
information taken from outside $S_1 \cup D_1$, the point process in
the remaining region (the new \emph{live zone}) is still a Poisson
process with unchanged intensity given the information within $S_1
\cup D_1$ --- as explained in Section \ref{DJBSim}.

We now grow a trapezium whose left--hand side located at $x=x_1$ has
length $y=y_1$. The trapezium expands until its right--hand side
first hits a seed $s_2$ (in the new live zone). The stopping set
formed is called $S_2$. It has an exponentially distributed area
$E_2$.

We proceed in this way, forming a sequence of stopping sets
(illustrated in Figure \ref{trapSS}) which, unlike those in Section
\ref{DJBSim}, do not form a nest. They do, however, have independent
exponentially--distributed areas and are part of a recursive
structure which we can exploit. It is also important to note that
the first $V$--type seed will provide the ray which blocks the test
ray. Without our introduction of dead zones, a complicated algorithm
rather like \texttt{block} would be required to check if a $V$--type
ray actually reaches the path of the test ray.

{\sc Remark 2:} \emph{Why is it that no seed in $D_1$ can influence
the distance $L$ travelled by the test seed; dead zone $V$--types
will either be blocked by an east--growing ray within the live zone
or, if they are not blocked, the test ray must have been intersected
at an earlier point. Dead zone $H$--types can never be in a
geometrical position to block live zone $V$--types. The same line of
reasoning applies to dead zones $D_2, D_3, ...$.}

{\sc Remark 3:} \emph{In Section 3 we mentioned that extraneous
randomisations, those not solely dependent on the Poisson point
process, might invalidate the key results from the stopping--set
theory. There is no such problem here with the stopping set $S_k$
itself, but we note that $D_k,\ k\geq1,$ depends for its existence
on an extraneous random feature --- namely the $H$ or $V$ mark of
seed $s_k$. This does not invalidate our comment above that the
point process in the current live zone outside $S_k\cup D_k$ is
unaffected by the information in $S_k\cup D_k$. For one thing the
seed marks are independent of each other and of the point process.
Furthermore, we only stop constructing dead zones when we have no
further need to observe the process at all. So the extraneous random
feature is not operative in our analysis.}

\vspace{.4cm}\textbf{The recursive structure commencing with a
generic live zone:} Suppose that we begin observing the process when
the live zone has left boundary of height $y$ and when we are about
to construct $S_n$. In Figure \ref{trapSS}, we draw the case $n=2$.
The probability density function for the length, $r$, of $S_n$'s
base, conditional on the height $y$ of its left boundary, follows
from the exponential distribution of $S_n$'s area $E_n$. It is
therefore:
\begin{equation*}
f(r\mid y) = \lambda (r+y)e^{-\frac{\lambda}{2}(r^2+2ry)}.
\end{equation*}
If the stopping seed $s_n$ for set $S_n$ is $V$--type, then its
south ray will be the first to intersect the test ray and the
process ends. Otherwise, another dead zone is created and further
trapezoidal stopping sets are formed until a $V$--type is met.

Let $X$ be the random variable equal to the horizontal distance
covered by stopping sets until the process comes to an end. The
density function of $X$, conditional on $y$ will be:
\begin{equation}
\label{conDensity}
g(x|y) = \frac{\lambda}{2}\left[ (x+y)e^{-\frac{\lambda}{2}(x^2+2xy)} + \int_0^\infty e^{-\frac{\lambda}{2}(r^2+2ry)}\left(\int_0^{r+y} g(x-r|u) du\right) dr\right],
\end{equation}
where $g(x|y)=0$ if $x<0$. The first term in the square brackets
accounts for the case where the first seed is V--type, and the
second term for the case where it is H--type and the process is
effectively re--started with a different boundary condition having
already covered some horizontal distance. We have taken
$q=\tfrac{1}{2}$, but the analysis can be carried out for general
$q$, producing a more complicated result. Note that the ray length
probability density function is $g(x\mid 0)$.

We define the moments of the conditional density:
\begin{equation*}
\mu_n(y) = \int_0^\infty x^n g(x|y) dx.
\end{equation*}
As mentioned before we will here compute $\mathbb{E}(L)=\mu_1(0)$ and $\mathbb{E}(L^2)=\mu_2(0)$, which from equation (\ref{conDensity}), satisfy:
\begin{align}
\label{mu1}
\mu_1(0) & = \sqrt{\frac{\pi}{2\lambda}}\left[1+\frac{\lambda}{2}\int_0^\infty \textrm{erfc}\left(\sqrt{\frac{\lambda}{2}}u\right)\mu_1(u) du\right] \\
\label{mu2}
\mu_2(0) & = \frac{2}{\lambda} + \frac{\lambda}{2}\sqrt{\frac{\pi}{2\lambda}}\int_0^\infty \textrm{erfc}\left(\sqrt{\frac{\lambda}{2}}u\right)\mu_2(u) du + \int_0^\infty e^{-\frac{\lambda}{2}u^2}\mu_1(u)du.
\end{align}
Our strategy is to find $\mu_1(y)$ and $\mu_2(y)$ up to an arbitrary constant, and then to determine the constant using equations (\ref{mu1}) and (\ref{mu2}). The first part of this process is most easily achieved by making use of the moment generating function:
\begin{equation*}
M_t(y)=\int_0^\infty e^{tx} g(x|y) dx,
\end{equation*}
which, from equation (\ref{conDensity}) satisfies
\begin{multline*}
M_t(y) =  \frac{1}{2} + \frac{1}{2}\sqrt{\frac{\pi}{2\lambda}} e^{\frac{(\lambda y-t)^2}{2\lambda}} \left\{ \textrm{erfc}\left(\frac{\lambda y-t}{\sqrt{2 \lambda}}\right)\left[t+\lambda \int_0^y M_t(u) du\right] \right. \\
\left. +\lambda \int_y^\infty  \textrm{erfc}\left(\frac{\lambda u-t}{\sqrt{2 \lambda}}\right) M_t(u) du \right\}.
\end{multline*}
This integral equation may be reduced to the differential equation:
\begin{equation*}
\frac{d^2M_t}{dy^2}-(\lambda y-t)\frac{dM_t}{dy} - \frac{\lambda}{2} M_t = - \frac{\lambda}{2}.
\end{equation*}
Expressing the left hand side as a series in $t$, and collecting coefficients of $t$ and $t^2$ we obtain differential equations satisfied by $\mu_1(y)$ and $\mu_2(y)$:
\begin{align}
\label{mu1Diff}
\mu_1''(y) - \lambda y \mu_1'(y) - \frac{\lambda}{2} \mu_1(y) & = 0 \\
\label{mu2Diff}
\mu_2''(y) - \lambda y \mu_2'(y) - \frac{\lambda}{2} \mu_2(y) & = -2 \mu_1'(y).
\end{align}
Clearly we must solve for $\mu_1(y)$ first.

\vspace{.4cm}\textbf{The first conditional moment:} Making the
change of variable $z=(\tfrac{\lambda}{2})^{\frac{1}{2}}y$ in
equation  (\ref{mu1Diff}) we obtain:
\begin{equation*}
\frac{d^2 \mu_1}{dz^2}-2z\frac{d \mu_1}{dz} - \mu_1=0.
\end{equation*}
If the coefficient of $\mu_1$ were a positive multiple of two, this would be Hermite's equation, solved by Hermite polynomials. Since this is not the case, we seek a series solution \cite{math73}:
\begin{equation*}
\mu_1(y(z))=\sum_{n=0}^\infty a_n z^n
\end{equation*}
and obtain the recurrence relation:
\begin{equation*}
a_{n+2}=\frac{2n+1}{(n+1)(n+2)}a_n.
\end{equation*}
This leads to the general solution:
\begin{equation*}
\mu_1(y(z)) = a_0 M(\tfrac{1}{4},\tfrac{1}{2},z^2) + a_1 z\ M(\tfrac{3}{4},\tfrac{3}{2},z^2),
\end{equation*}
where $M$ is Kummer's Function \cite{abram70}:
\begin{equation*}
M(a,b,z) = \sum_{n=0}^\infty \frac{(a)_n z^n}{(b)_n n!}.
\end{equation*}
Here we have used the Pochhammer symbol, defined by:
\begin{equation*}
(a)_n=a(a+1)(a+2)...(a+n-1),\ (a)_0=1.
\end{equation*}
The Kummer's functions diverge as $z \rightarrow \infty$, but we know that $\mu_1(y(z)) \rightarrow 0$ in that limit. This apparent paradox is resolved by noting that the two independent parts of the solution may be combined to form a Kummer's function of the second kind \cite{abram70}, defined by:
\begin{equation*}
U(a,b,z) = \frac{\pi}{\sin \pi b} \left[ \frac{M(a,b,z)}{\Gamma(1+a-b)\Gamma(b)} -
z^{1-b}\frac{M(1+a-b,2-b,z)}{\Gamma(a)\Gamma(2-b)}\right]
\end{equation*}
which tends to zero as $z \rightarrow \infty$. In terms of this function, the general solution is:
\begin{equation*}
\mu_1(y(z)) = A \times M(\tfrac{1}{4},\tfrac{1}{2},z^2) + B \times U(\tfrac{1}{4},\tfrac{1}{2},z^2).
\end{equation*}
It must be the case that $A=0$ in order to capture the right asymptotic behavior so, restoring the original variable $y$, the conditional moment must have the form:
\begin{equation}
\label{mu1Gen}
\mu_1(y) = B \times U(\tfrac{1}{4},\tfrac{1}{2},\tfrac{\lambda}{2}y^2).
\end{equation}
It now remains to compute $B$. We do this by substituting (\ref{mu1Gen}) into equation (\ref{mu1}). Making use of the result:
\begin{equation*}
\int_0^\infty \textrm{erfc}(u) U(\tfrac{1}{4},\tfrac{1}{2},u^2) du = \frac{\sqrt{2}}{\pi}[\Gamma(\tfrac{1}{4})-\sqrt{\pi}\Gamma(\tfrac{3}{4})],
\end{equation*}
together with $\Gamma(\tfrac{1}{4})\Gamma(\tfrac{3}{4})=\sqrt{2} \pi$ we find that:
\begin{equation*}
B=\frac{\sqrt{\pi}}{\sqrt{\lambda}\Gamma(\tfrac{3}{4})}.
\end{equation*}
We have now found $\mu_1(y)$, which gives us a compact analytic expression for the expected ray length:
\begin{align*}
\mathbb{E}(L) & = \mu_1(0)\\
& = \frac{\sqrt{\pi}}{\sqrt{\lambda}\Gamma(\tfrac{3}{4})} U(\tfrac{1}{4},\tfrac{1}{2},0) \\
& = \frac{\pi}{\sqrt{\lambda} \left(\Gamma(\tfrac{3}{4})\right)^2}\\
& \approx \frac{2.0920992}{\sqrt{\lambda}}
\end{align*}
For comparison, using the first 200 coefficients from Cowan and Ma's recurrence we obtain $\mathbb{E}(L) \approx 2.0920987$ when $\lambda=1$. As we discovered earlier, when $\lambda = 2$, the half model provides an approximation to the full model, having similar but simplified blocking effects and identical mean field behaviour. For this choice of $\lambda$ we obtain the exact half model result $\mathbb{E}(L) = 1.479337560$ to 7 decimal places, which differs from the accurate full model result (1.467535) by $0.7\%$. Compared with the mean field prediction: $\mathbb{E}(L) \approx \sqrt{2}$, which differs from the full model by $3.6\%$ this is a much closer approximation.

\vspace{.4cm}\textbf{The second conditional moment:} As for the
calculation of $\mu_1$, we make the change of variable
$z=(\tfrac{\lambda}{2})^{\frac{1}{2}}y$, but this time in equation
(\ref{mu2Diff}), obtaining
\begin{equation*}
\frac{d^2 \mu_1}{dz^2}-2z\frac{d \mu_1}{dz} - \mu_1= \frac{\sqrt{2\pi} z }{\lambda \Gamma(\tfrac{3}{4})} U(\tfrac{5}{4},\tfrac{3}{2},z^2),
\end{equation*}
where we have used the differential property \cite{abram70} : $U'(a,b,z)=-a\ U(a+1,b+1,z)$. We know the homogeneous part of the general solution to (\ref{mu2Diff}), so it remains to find a particular solution. We do this using variation of parameters, and begin by making the definitions:
\begin{align*}
f_1(z) & =  M(\tfrac{1}{4},\tfrac{1}{2},z^2) \\
f_2(z) & = U(\tfrac{1}{4},\tfrac{1}{2},z^2).
\end{align*}
The function $M$ has the differential property \cite{abram70}: $M'(a,b,z)=\tfrac{a}{b}M(a+1,b+1,z)$ which allows us to compute the Wronskian:
\begin{align*}
W(z) & = f_1(z)f_2'(z)-f_2(z)f_1'(z) \\
& = -\frac{z}{2}\left[M(\tfrac{1}{4},\tfrac{1}{2},z^2) U(\tfrac{5}{4},\tfrac{3}{2},z^2) + 2U(\tfrac{1}{4},\tfrac{1}{2},z^2)M(\tfrac{5}{4},\tfrac{3}{2},z^2)\right].
\end{align*}
We now define
\begin{equation*}
G(z,t)  = \frac{f_2(z)f_1(t)-f_1(z)f_2(t)}{W(t)}
\end{equation*}
in terms of which the particular integral is:
\begin{equation*}
f_p(z) = -\frac{\sqrt{2 \pi}}{\lambda\Gamma(\tfrac{3}{4})} \int_z^\infty G(z,t) t U(\tfrac{5}{4},\tfrac{3}{2},t^2) dt.
\end{equation*}
Discarding the divergent part of the solution, and restoring $y$, we have that:
\begin{equation*}
\mu_2(y) = C \times  U(\tfrac{1}{4},\tfrac{1}{2},\tfrac{\lambda}{2}y^2) -\frac{\sqrt{2\pi}}{\lambda \Gamma(\tfrac{3}{4})} \int_{z(y)}^\infty G(z(y),t) t U(\tfrac{5}{4},\tfrac{3}{2},t^2) dt,
\end{equation*}
where $C$ is an as yet undetermined constant. We find it by substituting our expression for $\mu_2(y)$ into equation (\ref{mu2}). Making use of the numerical integral:
\begin{align*}
K & = - \int_0^\infty \textrm{erfc}(z) \left[\int_{z}^\infty G(z,t) t U(\tfrac{5}{4},\tfrac{3}{2},t^2) dt \right] dz \\
& = 0.343146
\end{align*}
we find that:
\begin{equation*}
C=\frac{1}{\Gamma\left(\tfrac{3}{4}\right) \lambda}\left(\frac{\pi K}{\Gamma\left(\tfrac{3}{4}\right)}
+ 2\sqrt{2}\right).
\end{equation*}
Noting also that $f_p(0)=\tfrac{2}{\lambda}$ we have the final result that:
\begin{align*}
\mathbb{E}(L^2) & = \mu_2(0) \\
& = \frac{1}{\Gamma\left(\tfrac{3}{4}\right) \lambda}\left(\frac{\pi K}{\Gamma\left(\tfrac{3}{4}\right)}
+ 2\sqrt{2}\right) U(\tfrac{1}{4},\tfrac{1}{2},0) + f_p(0) \\
& = \frac{\pi ^{3/2} K+2 \Gamma \left(\frac{3}{4}\right)
   \left(\sqrt{2 \pi }+\Gamma \left(\frac{3}{4}\right)^2\right)}{\lambda
    \Gamma \left(\frac{3}{4}\right)^3}\\
& \approx \frac{6.37688}{\lambda}
\end{align*}
For comparison, using the first 200 coefficients from the Cowan--Ma
recurrence we obtain $\mathbb{E}(L) \approx 6.37686$ when
$\lambda=1$.

\section{Concluding comment}

Gilbert's tessellation is notoriously difficult to analyse, and even
the rectangular version studied by Mackisack and Miles remains
entirely without analytical results. In this paper we have shown
that the simplified rectangular model of Cowan and Ma, with only
half of the blocking rules of the Mackisack and Miles model, has a
number of tractable properties. As such, it is the only
Gilbert--style model, we believe, which has yielded any analytic
results.

\section*{Appendix: Expected length in the half model when $q \neq \tfrac{1}{2}$}

If $q$ is the proportion of seeds growing horizontally in the half model, then equation (\ref{conDensity}) becomes:
\begin{equation*}
g(x|y) = (1-q) \lambda (x+y)e^{-\frac{\lambda}{2}(x^2+2xy)} + q\lambda \int_0^\infty e^{-\frac{\lambda}{2}(r^2+2ry)}\left[\int_0^{r+y} g(x-r|u) du \right] dr.
\end{equation*}

\begin{figure}[h]
\psfrag{EL}{$\mathbb{E}(L_H)$} \psfrag{q}{$q$}
\begin{center}
\includegraphics[width=11 cm]{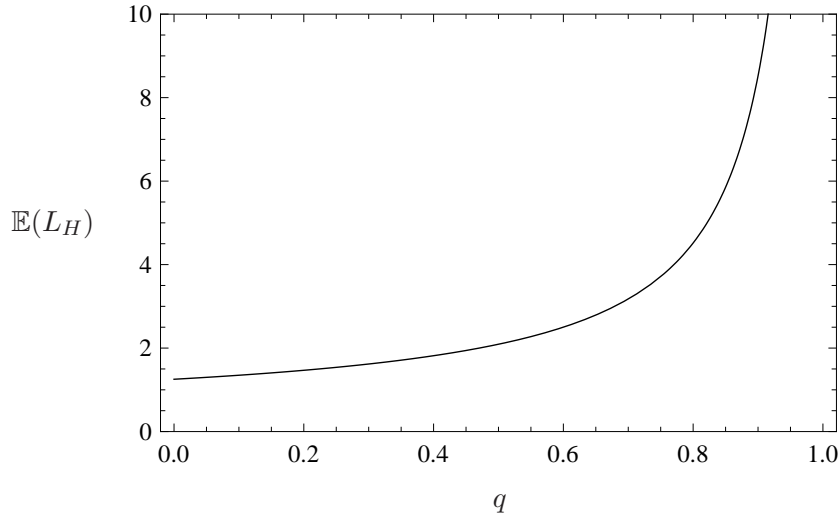}
\caption{\scriptsize \label{ELq} Expected horizontal length in the
half model as a function of $q$, the proportion of horizontal rays.
The seed density is $\lambda=1$.}
\end{center}
\end{figure}

The first moment of $g(x|0)$ may be found by similar methods to
those employed in the $q=\tfrac{1}{2}$ case. The expected length of
a horizontal ray is found to be:
\begin{equation*}
\mathbb{E}(L_H) = \sqrt{\frac{\pi}{\lambda}}\left[\sqrt{2}-\frac{q \Gamma\left(1-\tfrac{q}{2}\right)G_{3,3}^{2,3}\left(1\left|
\begin{array}{c}
 0,\frac{1}{2},\frac{q+1}{2} \\
 0,\frac{1}{2},-\frac{1}{2}
\end{array}
\right.\right)}{2^{q+\frac{1}{2}}\pi \Gamma(1-q)}\right]^{-1}
\end{equation*}
where $G$ is Meijer's G--Function \cite{grad07}. Figure \ref{ELq}
illustrates the function $\mathbb{E}[L_H]$ for the case $\lambda=1$.

\end{document}